\documentclass[fleqn]{article}
\usepackage{amssymb,latexsym,theorem}

\newcommand{\B}{\mathbf{B}}
\newcommand{\K}{\mathbb{K}}
\newcommand{\V}{\mathbb{V}}
\newcommand{\X}{\mathbb{X}}
\newcommand{\e}{\varepsilon}

\newcommand{\eps}{\varepsilon}
\newcommand{\di}{\mathrm{dim}}
\newcommand{\codi}{\mathrm{codim}}
\newcommand{\PG}{\mathrm{PG}}

\newcommand{\F}{\mathbb{F}}

\newcommand{\cP}{\mathcal{P}}

\newcommand{\spin}{{\rm spin}}
\newcommand{\ch}{{\rm char}}

\newtheorem{lemma}{Lemma}[section]
\newtheorem{defin}[lemma]{Definition}
\newtheorem{prop}[lemma]{Proposition}
\newtheorem{theo}[lemma]{Theorem}
\newtheorem{co}[lemma]{Corollary}

\def\N{\mathbb{N}}

\def\Aut{\mathrm{Aut}}
\def\PGL{\mathrm{PGL}}
\def\PGGL{\mathrm{P}\Gamma\mathrm{L}}
\def\DW{\mathrm{DW}}
\def\DQ{\mathrm{DQ}}
\def\dH{\mathrm{DH}}

\def\W{\mathrm{W}}
\def\Q{\mathrm{Q}}

\def\char{\mathrm{char}}

\def\gr{{\mbox{\rm{\footnotesize{gr}}}}}
\def\spin{{\mbox{\rm{\footnotesize{spin}}}}}

\def\d{\mathrm{d}}

\title{An outline of polar spaces: basics and advances}
\author{Ilaria Cardinali}
\date{ }

\begin{document}
\maketitle
\begin{abstract}
This paper is an extended version of a series of lectures on polar spaces given during the workshop and conference \lq Groups and Geometries\rq\, held at the Indian Statistical Institute in Bangalore in December 2012.
The aim of this paper is to give an overview of the theory of polar spaces focusing on some research topics related to polar spaces. We survey the fundamental results about polar spaces starting from classical polar spaces. Then we introduce and report on the state of the art on the following research topics: polar spaces of infinite rank, embedding polar spaces in groups and projective embeddings of dual polar spaces.
\end{abstract}

\section{Introduction}\label{Introduction}

Polar spaces have been widely studied over the years and many fundamental questions about polar spaces such as their classification and the description of remarkable examples have already been answered. The content of this paper has been the subject of three lectures given during the workshop and conference \lq Groups and Geometries\rq\, held at the Indian Statistical Institute in Bangalore in December 2012. This paper wants to give an overview of the theory of polar spaces starting from the classical polar spaces algebraically described by sesquilinear and pseudoquadratic forms and then enlarging the perspective to abstract polar spaces in which the classical polar spaces naturally sit. The aim of this general background is, naturally, the Classification Theorem of thick polar spaces of rank at least $3$.

Turning to research topics related to polar spaces, I decided to select three of them and report on  the state of the art in a way as complete as possible. Still, I have to say that this exposition is far of being a comprehensive treatment of all the possible branches of research related to polar spaces, but having to choose what to present due to space and time restriction, I focused on the material I found more appealing to me. So, the research topics I will deal with are the following: polar spaces of infinite rank; embeddings of polar spaces in groups; projective embeddings of dual polar spaces.

The organization of the paper is as follows.

In Section~\ref{basics} we provide the background of the theory of polar spaces. We start from classical polar spaces, hence introducing the algebraic notion of sesquilinear forms and pseudoquadratic forms (Section~\ref{sesquilinear forms} and Section~\ref{pseudoquadratic}). Then we define abstract polar spaces (Section~\ref{abstract polar spaces}) adopting the \lq 1-all axiom\rq\, and we conclude Section~\ref{basics} by stating the classification theorems in Section~\ref{fundamental theorems}. 

In Section~\ref{infinite rank} we deal with polar spaces of infinite rank. Several properties which always hold for polar spaces of finite rank fail to hold when we allow the rank of the polar space to be infinite. We report on the state of the art. Among other things, we describe an example of polar space of infinite rank. 


Section~\ref{sec4} focuses on embeddings of polar spaces in groups. After having provided the basic definitions of morphism and hull of an embedding, we report in Section~\ref{sec4-2} on two results which completely solve the problem of the non-abelian hull of a projective embedding of a polar space. 

Section~\ref{dual polar spaces} is about dual polar spaces and their projective embeddings. We give the definition of homogeneous, polarized, universal and minimal embedding of a dual polar space and then we give an update survey of the results on these embeddings for dual thick classical polar spaces. 

\section{Basic results}\label{basics}
The references books for the material of this section are Tits~\cite{Tits} and  Buekenhout-Cohen~\cite[Chapters 7-10]{BuekCohen}.

We start by introducing in Section \ref{sesquilinear forms} and in Section \ref{pseudoquadratic} the notion of sesquilinear and pseudoquadratic forms. In Section~\ref{abstract polar spaces} we deal with abstract polar spaces and the main theorems are reported in Section~\ref{fundamental theorems}.

\subsection{Sesquilinear forms}\label{sesquilinear forms}
\begin{defin}\label{def sesquilinear form}
Let $\K$ be a division ring and suppose that there exists a pair $( \sigma , \varepsilon )$ where $\sigma$ is an anti-automorphism of $\K$ and $\varepsilon \in \K^*$ such that


\begin{enumerate}
\item[(i)] $\varepsilon ^{\sigma}=\eps^{-1};$
\item[(ii)]$t^{{\sigma}^2}=\eps t \eps^{-1}.$
\end{enumerate}

Let $\V$ be a right vector space over $\K$. Then a function $\phi\colon \V\times \V \rightarrow \K$ is a {\it reflexive $(\sigma, \eps)$-sesquilinear form} (also {\it $(\sigma, \eps)$-sesquilinear form} for short) if

\begin{enumerate}
\item[(iii)] $\phi(x, y\alpha +z\beta)=\phi(x,y)\alpha +\phi(x, z)\beta \,\,\,\,\forall \alpha, \beta\in \K, \forall x,y\in \V$;
\item[(iv)] $\phi(y,x)=\phi(x,y)^{\sigma} \varepsilon , \,\,\,\,\forall x,y\in \V.$\\
\end{enumerate}
\end{defin}

The following are straightforward consequences of Definition~\ref{def sesquilinear form}.

\begin{enumerate}
\item[$(v)$] $\phi(\sum_i x_i \alpha_i,\sum_j y_j \beta_j )=\sum_{i,j}\alpha_i^{\sigma} \phi(x_i,y_j)\beta_j,\,\,\,\,\,\forall \alpha_i,\beta_j\in \K,\,\,\forall x_i,y_j\in \V$;
\item[$(vi)$] $\phi(x,x)=\phi(x,x)^{\sigma}\varepsilon \,\,\,\,\,\,\forall x\in \V.$
\end{enumerate}

We will adopt the following terminology and notation. If $a,b\in \V$ we will write $a\perp b$ when $\phi(a,b)=0$. We will say that $a$ is an {\it isotropic vector} if $a\perp a$ and we put $a^{\perp}:=\{x\in \V\colon x\perp a\}$. It is a direct consequence of $(iii)$ and $(iv)$ in Definition~\ref{def sesquilinear form} that $a^{\perp}$ is a linear subspace of $\V$.

If $X\subseteq \V$, then we will say that $X$ is {\it totally isotropic} if $\phi(x,y)=0 $ for every $x,y\in X$ and we define $X^{\perp}:=\cap_{x\in X}x^{\perp}.$ We can easily prove that $X^{\perp}=\langle X\rangle ^{\perp}$ and that $X^{\perp}$ is a linear subspace of $\V$.
Clearly,  $X$ is totally isotropic if and only if $X\subseteq X^{\perp}.$

 The following lemma states some elementary geometric properties which are immediate consequences of the axioms.

 \begin{lemma}\label{lem1}
 Let $a,b,c\in \V$ and $X\subseteq \V$. Then the following hold
 \item[(i)] $a\perp a$, $b\perp b$, $a\perp b\Rightarrow \langle a,b\rangle \subseteq \langle a,b\rangle^{\perp}$.
 \item[(ii)] $\langle X\rangle\subseteq \langle X\rangle^{\perp}\Leftrightarrow X\subseteq X^{\perp}.$
\item[(iii)]  $\langle b\rangle \not=\langle c\rangle$ and $\langle b,c\rangle\not\subseteq  a^{\perp}\Rightarrow \di(a^{\perp}\cap \langle b,c\rangle)=1.$
\item[(iv)] $a^{\perp}\not= \V \Rightarrow a^{\perp}$ is a hyperplane of $\V.$
 \end{lemma}

By Zorn's Lemma, every totally isotropic subspace is contained in a maximal totally isotropic subspace.
\begin{defin}\label{defin Witt index}
We will say that a $(\sigma,\e)$-sesquilinear form {\it admits a Witt index} if all maximal totally isotropic subspaces have the same dimension.
\end{defin}

Note that there are examples of $(\sigma,\e)$-sesquilinear forms which do not admit a Witt index as we will see in Section~\ref{infinite rank}.

 \begin{theo}\label{theo1}
 Let $\phi$ be a $(\sigma,\e)$-sesquilinear form.
If there exists a maximal isotropic subspace of finite dimension then $\phi$ admits a Witt index.
 \end{theo}

\begin{defin}
The subspace $\V^{\perp}=\{a\colon a^{\perp}=\V\}$ is called the {\em Radical of $\phi$} and denoted by $Rad(\phi).$  The form $\phi$ is degenerate (respectively, non-degenerate) if $Rad(\phi)\not= \{0\}$ (respectively, $Rad(\phi)= \{0\}$).
\end{defin}

Note that $Rad(\phi)\subseteq S$ for every maximal totaly isotropic subspace $S$ of $\V.$

The following proposition shows that we may always restrict ourselves to dealing only with non-degenerate sesquilinear forms. 

\begin{prop}\label{prop2}
Suppose that $\phi$ is a degenerate $(\sigma, \eps)$-sesquilinear form and let $W$ be a complement of $Rad(\phi)$ in $\V$. Then the form $\bar{\phi}\colon W\times W \rightarrow \K,\,\bar{\phi}:=\phi|_{W\times W}$ induced by $\phi$ on $W$ is non-degenerate and the function mapping every totally isotropic subspace $S$ of $\V$ containing $Rad(\phi)$ to $S\cap W$ is an isomorphism between the partially ordered set of the totally isotropic subspaces of $\V$ containing $Rad(\phi)$ to
the partially ordered set of the totally isotropic subspaces contained in  $W.$
 \end{prop}

If $\phi$ is a non-degenerate $(\sigma, \eps)$-sesquilinear form and $\X$ is a totally isotropic subspace then it is possible to prove that the codimension of $\X^{\perp}$ in $\V$ equals the dimension of $\X.$ Moreover, if $\X$ is finitely generated then $\X^{\perp \perp}=\X.$

An important class of sesquilinear forms is given by the so-called {\it trace-valued} sesquilinear forms.
\begin{defin}
A  $(\sigma, \eps)$-sesquilinear form $\phi$ is trace-valued if $\phi(x,x)\in \{t+t^{\sigma}\eps \}_{t\in \K}$ for every $x\in \V.$
\end{defin}

 \begin{lemma}\label{lem3}
If $\phi$ is  a trace-valued $(\sigma, \eps)$-sesquilinear form  and $a\not\in Rad(\phi)$ is an isotropic point for $\phi$ then for every line $l$ through $a$ such that $l\not\subseteq a^{\perp}$ there exists an isotropic point $b\not=a$ such that $b\in l.$
\end{lemma}
The following result is of fundamental importance.
 \begin{theo}\label{thm4}
Let $\phi$ be a $(\sigma, \eps)$-sesquilinear form  different from the null form and suppose that there exist isotropic points. Then $\V$ is spanned by the isotropic points if and only if both the following hold\\
(a) $\phi$ is trace-valued;\\
(b) there exist isotropic points not contained in the Radical of $\phi.$
\end{theo}

We will now define two subsets $\K_1$ and $\K_2$ of $\K$ as
 \begin{equation}
\K_1:=\{t\colon t=t^{\sigma}\eps\}\,\,\,{\mathrm{and}}\,\,\,\,\,\K_2:=\{t+t^{\sigma}\eps \}_{t\in \K}.
 \end{equation}

We immediately verify that $\K_2\subseteq \K_1.$ By Definition~\ref{def sesquilinear form} we have that $\phi(x,x)\in \K_1$ for every $x\in \V$ and that $\phi$ is trace-valued if and only if 
$\phi(x,x)\in \K_2$ for every $x\in \V.$

Denote by $\ch(\K)$ the characteristic of $\K$ and by $Z(\K)$ the center of $\K.$ Then
\begin{prop}\label{prop5}
If (i) $\ch(\K)\not= 2$ or (ii) $\ch(\K)= 2$ and $\sigma|_{Z(\K)}\not=id$ then $\K_2=\K_1$, namely $\phi$ is trace-valued.
\end{prop}

By Proposition~\ref{prop5}, if $\K$ is a field the only  $(\sigma, \eps)$-sesquilinear forms which are not trace-valued are for $\ch (\K)=2$ and $\sigma=id$, which force $\eps=1$.

 Suppose now $\ch(\K)=2$  and $\sigma|_{Z(\K)}=id.$ Suppose also that the sesquilinear form $\phi$ is not trace-valued and that there exist some isotropic vectors. Note that $\phi$ is not the null form since the null form is always trace valued. Let $\V_0$ be the subspace of $\V$ spanned by the isotropic vectors. By Theorem~\ref{thm4}, $\V_0$ is a proper subspace of $\V.$ Let $\phi_0:=\phi|_{\V_0\times \V_0}$ be the restriction of $\phi$ to $\V_0\times \V_0.$  By Theorem~\ref{thm4}, if $\phi_0$ is not the null form, thus $\phi_0$ is trace-valued. Moreover, every totally isotropic subspace of $\V$ is contained in $\V_0$.

So, given a sesquilinear form $\phi$ which is not trace-4valued, it is always possible to consider the associated non-degenerate trace-valued sesquilinear form $\phi_0$ in the way described above.

For the rest of this section we will only deal with $(\sigma,\varepsilon)$-sesquilinear forms admitting a finite Witt index $n\in \mathbb{N},$ according to Definition~\ref{defin Witt index}.

\begin{lemma}\label{lem6}
If $A$ is a maximal totally isotropic subspace of $\V$ and $p\not\in A$ is an isotropic point then $\langle A\cap p^{\perp},p\rangle$ is a maximal totally isotropic subspace.
\end{lemma}

 \begin{theo}\label{thm7}
Suppose $\phi$ is a trace-valued sesquilinear form of finite Witt index. Then, for every maximal totally isotropic subspace $A$ of $\V$ there exists a maximal totally isotropic subspace $B$ such that $A\cap B=Rad(\phi).$
\end{theo}
By Theorem \ref{thm7}, we immediately have
 \begin{co}\label{co8}
Suppose $\phi$ is a trace-valued sesquilinear form of finite Witt index $n$. Then $2n\leq \di(\V).$
\end{co}
We say that $\phi$ has maximal Witt index if $\di(\V)=2n$ or $\di(\V)=2n+1.$

Given a $(\sigma, \eps)$-sesquilinear form $\phi$ defined over a division ring $\K$, we can consider the Gram matrix $M$ of $\phi$. We recall that if $(e_i)_{i\in I}$ is a basis of $\V$, then the entries $M_{ij}$ of $M$ are definied as $M_{ij}:=\phi(e_i,e_j)$ for every $i,j\in I.$ It is well known that every  $(\sigma, \eps)$-sesquilinear form is completely determined by its Gram matrix and by $\sigma.$

Theorem \ref{thm7} will be used to construct a basis $\B$ of $\V$ so that the Gram matrix of $\phi$ with respect to $\B$ has a canonical form. We briefly explain how to proceed.

Suppose the form $\phi$ is non-degenerate and trace-valued with finite Witt index $n.$ By Theorem \ref{thm7}, there exist two maximal totally isotropic subspaces $A$ and $B$ of $\V$ such that $A\cap B=0.$ Suppose $a_1,\dots, a_n$ is a basis of $A.$

\begin{lemma}\label{lem9}
For every $k=i,\dots, n$, then $B\cap (\cap_{i\not=k}a_i^{\perp})$ is a point not in $a_k^{\perp}.$
\end{lemma}

Let $b_k$ be a vector generating  $B\cap (\cap_{i\not=k}a_i^{\perp}).$ By Lemma \ref{lem9}, $\phi(a_k,b_k)\not= 0$ for every $k=1,\dots, n$ and $\phi(a_i,b_k)= 0$ for $i\not= k.$ It is always possible to rescale $b_k$ so that $\phi(a_k,b_k)= 1.$

\begin{lemma}\label{lem10}
The vectors $b_1,\dots, b_n$ form a basis of $B.$
\end{lemma}

Since $\V=(A+B)\oplus (A+B)^{\perp}$, we can construct the basis $\B$ by taking the vectors $a_1,\dots, a_n, b_1, \dots, b_n $ and then completing this set with a basis of $(A+B)^{\perp}.$ The Gram matrix of $\phi$ with respect to $\B$ is then
$\left(\begin{array}{lll}
0_n & I_n & 0\\
\eps I_n & 0_n & 0\\
0 & 0 & M_0\\
\end{array} \right)
$

where $0_n$ is the $(n\times n)$-null matrix, $I_n$ is the $(n\times n)$-identity matrix and $M_0$ is the matrix representing the restriction of $\phi$ to the space $(A+B)^{\perp}.$

\begin{lemma}\label{lem11}
$(A+B)^{\perp}$  does not contain any isotropic point.
\end{lemma}

Note that if the division ring is finite or countable then we can always choose an orthogonal basis of $(A+B)^{\perp}$ so that $M_0$ is a diagonal matrix.

\subsubsection{Examples of non-degenerate trace-valued sesquilinear forms}\label{esempi}
\textbf{Alternating forms.}
If $\ch (\K)\not=2$ then a trace-valued $(\sigma,\eps)$-sesquilinear form $\phi$  is {\it alternating} if $(\sigma,\eps)=(id,-1).$
If $\ch (\K)=2$ then a trace-valued $(\sigma,\eps)$-sesquilinear form $\phi$  is alternating if $(\sigma,\eps)=(id,1)$ and $\phi(x,x)=0$ for every $x\in \V.$

If $\phi$  is alternating, all the vectors in  $\V$ are isotropic. Non-degenerate alternating forms of Witt index $n$ exist only in vector spaces of dimension $2n.$ The canonical Gram matrix of a non-degenerate alternating form (up to rescaling the form) is  $\left(\begin{array}{ll}
0_n & I_n \\
-I_n & 0_n  \\
\end{array} \right).
$

\textbf{Symmetric forms.}
A $(\sigma,\eps)$-sesquilinear form $\phi$  is {\it symmetric} if $(\sigma,\eps)=(id,1).$
Suppose $\K$ is a field and that every element of $\K$ is a square. If $\ch(\K)\not=2$ then the non-degenerate symmetric bilinear forms with Witt index $n$ exist only if $\di(\V)=2n$ or $\di(\V)=2n+1.$
If $\ch(\K)=2$ then the unique non-degenerate trace-valued symmetric bilinear forms are alternating.
\medskip

Let $\ch(\K)\not=2.$ If the multiplicative subgroup of the square elements of $\K$ has index $2$ in the multiplicative group $\K^*$ of $\K$ then we have the following possibilities.

(a) $\di(\V)=2n$ and the canonical Gram matrix of a non-degenerate symmetric form (up to rescaling the form) is  $\left(\begin{array}{ll}
0_n & I_n \\
I_n & 0_n  \\
\end{array} \right).
$
\\

(b) $\di(\V)=2n+1$ and the canonical matrix of a non-degenerate symmetric form (up to rescaling the form) is  $\left(\begin{array}{lll}
0 & I_n & 0\\
I_n & 0 & 0 \\
0&0& 1\\
\end{array} \right).
$
\\

(c) $\di(\V)=2n+2$ and the canonical matrix of a non-degenerate symmetric form (up to rescaling the form) is  $\left(\begin{array}{llll}
0_n & I_n & 0 & 0\\
I_n & 0_n & 0& 0\\
0 & 0 & 1& 0\\
0 & 0 & 0& \eta\\
\end{array} \right)
$ where $\eta$ is a non-square of $\K$.
\\

(d) $\di(\V)\geq 2n+2$ and the canonical matrix of a non-degenerate symmetric form (up to rescaling the form) is  $\left(\begin{array}{lll}
0_n & I_n & 0 \\
I_n & 0_n & 0\\
0 & 0 & M_0\\
\end{array} \right)
$ where $-1$ is a non-square and $M_0$ is the identity matrix of rank $\dim(\V)-2n$.
\\

If $\K$ is a finite field of odd characteristic then only cases (a), (b), (c) occur.

Note that if $\sigma$ is the identity then $\eps=\pm 1$ hence the previously described alternating and symmetric forms are the unique possible examples of $(\sigma,\eps)$-sesquilinear forms.

\textbf{Hermitian and anti-hermitian Forms.}
A $(\sigma,\eps)$-sesquilinear form $\phi$  is {\it hermitian} if $\sigma\not=id$ and $\eps=1.$ A $(\sigma,\eps)$-sesquilinear form $\phi$  is {\it anti-hermitian} if $\sigma\not=id$ and $\eps=-1.$

As an example, we will consider the case in which $\K$ is a finite field $\F_q$ and $\phi$ is a  $(\sigma,\eps)$-sesquilinear form.  Then, $q=q_0^2$ and we can always suppose $\eps=1$ (see Theorem~\ref{theo12}). If $\phi$ is a non-degenerate trace-valued $(\sigma, 1)$-hermitian form with finite Witt index $n$ then the following two possibilities might occur:\\

(i) $\di(\V)=2n$ and the canonical matrix of $\phi$ (up to rescaling the form) is  $\left(\begin{array}{ll}
0_n & I_n \\
I_n & 0_n  \\
\end{array} \right)$;
\\

(ii) $\di(\V)=2n+1$ and the canonical matrix of $\phi$ (up to rescaling the form) is  $\left(\begin{array}{lll}
0_n & I_n & 0\\
I_n & 0_n & 0 \\
0 & 0 & 1 \\
\end{array} \right).$
\\

The following theorem classifies all reflexive trace-valued $(\sigma, \eps)$-sesquilinear forms. Recall that two sesquilinear forms $\phi$ and $\psi$ are proportional if there exists  $\lambda\in \K^*$ such that $\phi=\lambda \psi.$

\begin{theo}\label{theo12}
If $\K$ is a division ring  then every reflexive trace-valued $(\sigma, \eps)$-sesquilinear form is either symmetric, alternating or proportional to a hermitian form.
\end{theo}

\subsection{Pseudoquadratic forms}\label{pseudoquadratic}
In this Section we introduce pseudoquadratic forms.
The motivation to introduce them is due to the fact that there exists an important class of classical polar spaces which can not be defined by means of trace-valued sesquilinear forms.

As in Section \ref{sesquilinear forms}, denote by $\K$, $\sigma$ respectively $\eps$ a division ring, an anti-automorphism of $\K$ and an element of $\K^*$  satisfying $(i)$ and $(ii)$ of Definition~\ref{def sesquilinear form}.
By Proposition \ref{prop5}, if $\K$ is a division ring of characteristic different from $2$ or $\ch(\K)=2$ and   $\sigma|_{Z(\K)}=id$ then a $(\sigma,\eps)$-sesquilinear form is trace-valued.

Put
\begin{equation}
\K_{\sigma, \eps}:=\{ t-t^{\sigma}\eps \}_{t\in \K}.
\end{equation}
Then $\K_{\sigma, \eps}$ is a subgroup of the additive group of  $\K.$ 
\begin{theo}\label{thm13}
$\K_{\sigma, \eps}=\K$ if and only if $\ch(\K)\not= 2$ and $(\sigma, \eps)=(id, -1).$
\end{theo}
Clearly, if $(\sigma, \eps)=(id,1)$ then $\K_{\sigma, \eps}=0.$ Moreover, for every  $\lambda\in \K$ and every $\tau\in \K_{\sigma, \eps}$ we have $\lambda^{\sigma}\tau\lambda \in \K_{\sigma, \eps}.$

\begin{defin}\label{pseudo quadratic form}
A function $f\colon \V\rightarrow \K/\K_{\sigma, \eps}$ is a $(\sigma, \eps)$-pseudoquadratic form (quadratic form, when $(\sigma, \eps)=(id,1)$) if
\item[(i)] $f(tx)=t^{\sigma}f(x)t\,\,\,\,\forall x\in \V,\,\,\forall t\in \K;$
\item[(ii)] $f(x+y)=f(x)+f(y)+(\phi(x,y)+\K_{\sigma, \eps})$ for every $x,y\in \V$ where $\phi$ is a suitable $(\sigma, \eps)$-sesquilinear form.
\end{defin}

By Definition \ref{pseudo quadratic form} it easily follows that the form $\phi$ is trace-valued and that $f(0)=0$ and $f(-x)=f(x)$ for every $x\in \V.$

\begin{theo}\label{theo13}
If $\K_{\sigma, \eps}\not=\K$ then the sesquilinear form $\phi$ is uniquely determined by the pseudoquadratic form $f.$
\end{theo}

If $\K_{\sigma, \eps}=\K$, i.e. 
if $\ch(\K)\not=2$ and $(\sigma,\e)=(id,-1)$ (hence we are dealing with alternating forms in odd characteristic) then the associated pseudoquadratic  form is null.

From hereon suppose $\K_{\sigma, \eps}\not=\K$. By Theorem \ref{theo13} the sesquilinear form $\phi$ is uniquely determined by the pseudoquadratic form $f$. The sesquilinear form $\phi$ associated to the pseudoquadratic form $f$ is called the {\it sesquilinearization of $f$}.
The following theorem states under which conditions a sesquilinear form uniquely defines a pseudoquadratic form.

\begin{theo}\label{theo14}
(a) If $\ch(\K)\not= 2$ then $f(x)=\phi(x,x)/2 + \K_{\sigma, \eps},\,\,\forall x\in \V$;\\
(b) If $\ch(\K)= 2$ and $\sigma|_{Z(\K)}\not=id$ then \[f(x)=\phi(x,x)/(1+(t^{\sigma}/t)^2) + \K_{\sigma, \eps},\,\,\,\,\forall x\in \V\] where $t\in Z(\K)$ with $t^{\sigma}\not= t.$
\end{theo}
\subsubsection{Singularity}
Suppose  $\K_{\sigma, \eps}\not=\K$ and that $f$ is a $(\sigma, \eps)$-pseudoquadratic form with $\phi$ its sesquilinearization.

\begin{defin}\label{singularity}
A point $p$ in $\V$ is singular if $f(p)=0.$ A subspace $S$ of $\V$ is totally singular if $f(x)=0$ for every $x\in S.$
\end{defin}

The following result enlightens the link between totally singular subspaces and totally isotropic subspaces.
\begin{theo}\label{theo15}
If $S$ is a totally singular subspace for $f$ then $S$ is a totally isotropic subspace for $\phi.$ In particular, all the points of  $\V$ which are singular for $f$ are isotropic for $\phi.$
\end{theo}

By Theorem \ref{theo14}, if $\ch(\K)\not= 2$ or  $\ch(\K)= 2$ and $\sigma|_{Z(\K)}\not=id$ then every isotropic point for $\phi$ is also singular for $f.$ Under these hypothesis, the theory of totally isotropic subspaces and the theory of totally singular subspaces are the same.

From hereon we assume $\ch(\K)= 2$ and $\sigma|_{Z(\K)}=id.$ 

Denote by $S(f)$ the set of singular points of $\V.$ Then for every totally isotropic subspace $S$ of $\V$ we can prove that $S\cap S(f)$ is a (possibly empty) subspace of $S.$ In particular, $S(f)\cap Rad(\phi)$ is a subspace of $Rad(\phi).$

\begin{defin}\label{singularity2}
 The defect of a pseudoquadratic form $f$ is the codimension of $S(f)\cap Rad(\phi)$ as a subspace of $Rad(\phi).$ If $S(f)\cap Rad(\phi)=0$ [respectively, $S(f)\cap Rad(\phi)\not=0$] then $f$ is non-singular [respectively, $f$ is singular].
\end{defin}

Many results proved in Section~\ref{sesquilinear forms} for (non-degenerate) sesquilinear forms have an analogue for (non-singular) quadratic forms.

We say that a $(\sigma,\e)$-pseudoquadratic form admits a Witt index if all maximal singular subspaces have the same dimension. Note that not every $(\sigma,\e)$-pseudoquadratic form admits a Witt index.

The following is the analogue of Theorem \ref{theo1}
\begin{theo}\label{theo16}
 Let $f$ be a $(\sigma,\e)$-pseudoquadratic form.
If there exists a maximal singular subspace for $f$ of finite dimension then $f$ admits a Witt index.
 \end{theo}

The following is the analogue of Theorem \ref{thm4}
\begin{theo}\label{theo17}
If $S(f)\not\subseteq Rad(\phi)$ then $\langle S(f)\rangle =\V.$
\end{theo}

Given a singular pseudoquadratic form it is always possible to consider the non-singular pseudoquadratic form associated to it. Indeed

\begin{prop}\label{prop18}
Suppose $f$ is a singular quadratic form and let $W$ be a complement of $S(f)\cap Rad(\phi)$ in $\V.$ Then the form $f|_W$ induced by $f$ on $W$ is non-singular, its sesquilinearization is $\phi |_{W\times W}$ and the function $\pi_W$ mapping every totally singular subspace $S$ of $\V$ containing  $S(f)\cap Rad(\phi)$ to $S\cap W$ is an isomorphism from the partially ordered set of the totally singular subspaces of $\V$ containing $S(f)\cap Rad(\phi)$ to the partially ordered set of the totally singular (for $f|_{W}$) subspaces of $W.$
\end{prop}
The following is the analogue of Theorem \ref{thm7}

\begin{theo}\label{thm19}
Suppose $S(f)\not\subseteq Rad (\phi)$ and that $f$ has finite Witt index. Then, for every maximal totally singular subspace $A$ of $\V$ there exists a maximal totally singular subspace $B$ such that $A\cap B=Rad(\phi)\cap S(f).$
\end{theo}

Denote by $n_f$ the Witt index of a quadratic form $f$ and by $n_{\phi}$ the Witt index of the sesquilinearization $\phi$ of  $f.$
Then \begin{co}\label{co20}
\[\di(Rad(\phi))+n_f\leq n_{\phi}+  \di (Rad(\phi)\cap S(f)).\]
\end{co}

To conclude this section we focus on a special example.
\\

\textbf{Suppose $\K$ be a perfect field of characteristic $2$.}
We want to classify all non-singular quadratic forms (hence $(\sigma, \eps)=(id, 1)$). 

\begin{theo}\label{thm20}
If $f$ is a non-singular quadratic form and $\phi$ is its sesquilinearization then $\di(Rad(\phi))\leq 1.$ If $\di(Rad(\phi))= 1$ (so, $Rad(\phi)$ is a projective point) then every line of $\PG(\V)$ through the point $Rad(\phi)$ contains only one singular point.
\end{theo}
The following two cases might occur\\

Case 1. \fbox{$Rad(\phi)=0.$} Then the maximal singular subspaces of $\V$ for $f$ are either maximal isotropic subspaces for $\phi$ or they are hyperplanes of maximal isotropic subspaces for $\phi.$ Put $n=n_{\phi}.$ Then $\di(\V)=2n$ and the sesquilinearization of $f$ is a non-degenerate alternating form.\\

Case 1.1. Hyperbolic case: $n_{\phi}=n_f=n.$ It is possible to choose a suitable basis $\B=\{a_1,\dots, a_n,b_1,\dots b_n\}$ of $\V$ so that the quadratic form $f$ has the following expression, where $x= \sum_{i=1}^{n} x_ia_i+\sum_{i=1}^{n} y_ib_i$ is an arbitrary vector of $\V$:
\[f(\sum_{i=1}^{n} x_ia_i+\sum_{i=1}^{n} y_ib_i)=\sum_{i=1}^{n} x_iy_i.\]

Case 1.2. Elliptic case: $n_f=n_{\phi}-1=n-1.$  It is possible to choose a suitable basis $B=\{a_1,\dots, a_n,b_1,\dots b_n\}$ of $\V$ so that the quadratic form $f$ has the following expression, where $x= \sum_{i=1}^{n} x_ia_i+\sum_{i=1}^{n} y_ib_i$ is an arbitrary vector of $\V$:
\[f(\sum_{i=1}^{n} x_ia_i+\sum_{i=1}^{n} y_ib_i)=\sum_{i=1}^{n-1} x_iy_i +x_n^2+\eta x_ny_n+y_n^2\]
where $\eta\in \K$ so that the equation $t^2+\eta t +1=0$ has no solution in $\K.$\\

Case 2. \fbox{$\di (Rad(\phi))=1.$} Then the maximal isotropic subspaces of $\V$ for $\phi$ are of the form $Rad(\phi)\oplus S$ where $S$ is a maximal totally singular subspace for $f.$ So, $n=n_{\phi}=n_f+1.$ The point $Rad(\phi)$ is said to be the {\it nucleus} of $f$ and it is characterized as the unique non-singular point such that every line through it meets $S(f)$ in exactly one point.

It is possible to choose a suitable basis $\B=\{a_1,\dots, a_n,b_1,\dots b_n,c\}$ of $\V$ so that the quadratic form $f$ has the following expression, where $x= \sum_{i=1}^{n} x_ia_i+\sum_{i=1}^{n} y_ib_i+zc$ is an arbitrary vector of $\V$:
\[f(\sum_{i=1}^{n} x_ia_i+\sum_{i=1}^{n} y_ib_i+zc)=\sum_{i=1}^{n} x_iy_i +z^2.\]





\subsection{Abstract polar spaces} \label{abstract polar spaces}
We recall that a partial linear space $(P,L)$ is a point-line geometry where $P$ is the pointset, $L$ is the lineset, any line is incident with at least two points and any pair of distinct points is incident with at most one line.

If $(P,L)$  is a partial linear space we will adopt the following notation. If $a,b\in P$ then $a\perp b$ means that $a$ and $b$ are collinear. Put $a^{\perp}:=\{p\colon p\perp a\}\cup \{a\}$ and if $X\subseteq P$ then $X^{\perp}:=\cap_{x\in X} x^{\perp}$.

\begin{defin}[Buekenhout-Shult~\cite{Buek-Shult}]\label{def abstract polar spaces}
A partial linear space $\cP=(P,L)$ is a non-degenerate ordinary polar space if
\item[(i)] $P^{\perp}=\emptyset.$
\item[(ii)] For each point $p\in P$ and each line $l\in L,$ the set of points of $l$ which are collinear with $p$ is either a singleton or $l.$
\item[(iii)] Each line has at least three points.
\end{defin}

Condition $(ii)$ is also called the \lq 1-all axiom\rq\, for polar spaces.

By Zorn's Lemma, every singular subspace is contained in a maximal singular subspace.

\begin{defin}
Let $\cP=(P,L)$ be a non-degenerate ordinary polar space and $S\subseteq P.$ Then $S$ is singular if $S\subseteq S^{\perp}.$\\
A polar space $\cP$ admits a Witt index if all maximal singular subspaces have the same dimension.

\end{defin}

\begin{theo}\label{thm21}
Let $\cP$ be a non-degenerate ordinary polar space. If there exists a maximal chain of singular subspaces of $\cP$ of finite length then every maximal chain of singular subspaces has finite length and all maximal chains have the same length.
\end{theo}

\begin{defin}
Let $\cP$ be a non-degenerate ordinary polar space. If $\cP$ admits a maximal chain of singular subspaces of finite length then, in view of Theorem \ref{thm21}, we define the rank of $\cP$ to be the length of a maximal chain of non-empty singular subspaces.

If $\cP$ does not admit any maximal chain of singular subspaces of finite length then we say that the rank of $\cP$ is infinite.
\end{defin}

A polar space of rank $2$ in which for every point $p$ and every line $l$ with $p$ not on $l$ there exists a unique point on $l$ collinear with $p$ is called {\it generalized quadrangle}.
At the other extreme, a projective space is an example of polar space where given a point $p$ and a line $l$ with $p$ not on $l$, all points on $l$ are collinear with $p.$

Throughout this section, $\cP$ is assumed to be a non-degenerate ordinary polar space of finite rank. From Definition~\ref{def abstract polar spaces}, we can deduce many geometric results some of which will be recalled in the following theorems. 

\begin{theo}\label{thm22}
Let $\cP=(P,L)$ be a non-degenerate ordinary polar space of finite rank. Then the following hold
\item[(i)] Every singular subspace of $\cP$ is a projective space.
\item[(ii)] If $S$ is a maximal subspace of $\cP$ and $p$ a point not in $S$ then $p^{\perp}\cap S$ is a hyperplane of $S$ and there exists a unique maximal subspace $S'=\langle p^{\perp}\cap S, p \rangle$ such that $S\cap S'=p^{\perp}\cap S.$ 
\item[(iii)] There exist two disjoint maximal subspaces.
\end{theo}

The following theorem describes an important family of polar spaces.
\begin{theo}\label{thm23}

(i) Let $\phi$ be a non-degenerate reflexive trace-valued $(\sigma,\eps)$- sesquilinear form (of finite Witt index) over a vector space $\V$. Then the sets of totally isotropic points and lines of $\PG(\V)$ with respect to $\phi$ is a non-degenerate ordinary polar space (of finite rank). \\

(ii) Let $f$ be a non-singular pseudoquadratic form (of finite Witt index) over a vector space $\V$. Then the sets of totally singular points and lines of $\PG(\V)$ with respect to $f$ is a non-degenerate ordinary polar space (of finite rank). \\
\end{theo}

\begin{defin}
A polar space which is obtained from a non-degenerate reflexive trace-valued $(\sigma,\eps)$- sesquilinear or from a non-singular pseudoquadratic form is called a classical polar space.
\end{defin}

According to Subsection \ref{esempi}, we will say that a classical polar space  is symplectic if it can be obtained from an alternating sesquilinear form; it is orthogonal if it can be obtained from a symmetric  sesquilinear form or from a non-singular quadratic form; it is hermitian if it can be obtained from a hermitian or anti-hermitian sesquilinear form.

The most natural question to ask now is whether there are some non-degenerate ordinary polar spaces of finite rank which are not classical. The affirmative answer is given in the next section.

\subsection{Fundamental  theorems}\label{fundamental theorems}
\begin{defin}
A polar space $\cP=(P,L)$ is embeddable if there exists a vector space $\V$ and an injective map $\xi\colon P\rightarrow \PG(\V)$ such that
\item[(E1)] $\xi(l)$ is a line of $\PG(\V)$ for every line $l\in L$;
\item[(E2)] $\langle \xi(P)\rangle=\PG(\V).$
\end{defin}

Note that every classical polar space is embeddable and an embeddable polar space is necessarily ordinary and all its planes are desarguesian.

\begin{theo}[Tits~\cite{Tits}, Buekenhout and Lef\`evre~\cite{BL}]\label{thm24} Every embeddable polar space of rank $n\geq 2$ is classical.
\end{theo}

\begin{defin}
An ordinary polar space of rank $n$ is thick if every singular subspace of dimension $n-2$ is contained in at least three singular subspaces  of dimension $n-1$ (i.e. maximal subspaces).
A polar space of rank $n$ is top-thin if every singular subspace of dimension $n-2$ is contained in exactly two singular subspaces  of dimension $n-1$ (i.e. maximal subspaces).
\end{defin}

\begin{theo}[Tits~\cite{Tits}]\label{thm25}
\item[(1)] Every ordinary polar space of rank $n\geq 4$ is embeddable (hence classical).
\item[(2)] A thick polar space of rank $3$ is embeddable if and only if  its planes are desarguesian.
\item[(3)] An ordinary  top-thin polar space of rank $3$ is embeddable if and only if its planes are Pappian. \\
Moreover
\item[(4)] There exists a unique family of non-embeddable thick polar spaces of rank $3.$ The planes of this polar space are Moufang planes.
\item[(5)] Every  ordinary, top-thin polar space of rank $3$ is obtained as the Grassmannian of lines of a projective space $\PG(3,\K).$
\end{theo}

A classical reference for a description of the non-embeddable polar space of rank 3 as in (4) of Theorem~\ref{thm25} is given in Tits~\cite[Chapter 9]{Tits}. See also M\"{u}hlerr~\cite{Muhlherr} for another more geometrical description  of it. Note that a new description of the family as in in (4) of Theorem~\ref{thm25} has been recently given by De Bruyn-Van Maldeghem~\cite{DB-VM}.
\\

By Wedderburn's Theorem, every finite division ring is commutative and by Artin and Zorn's Theorem, every finite projective plane of Moufang type is Pappian. So, we have

\begin{co}\label{thm26}
Every ordinary finite polar space of rank at least $3$ is classical.
\end{co}

\section{Polar spaces of infinite rank}\label{infinite rank}
In this section we will focus on polar spaces of infinite rank. Many properties of polar spaces of finite rank fail to hold for polar spaces of infinite rank.
We will first give an overview of the properties which do not hold anymore in a polar space of infinite rank and then we will provide in Subsection \ref{example-infinite rank} an example of polar space of infinite rank where some of these phenomena actually happen. We will then stress on the open problems regarding polar spaces of infinite rank.

The material of this section is mainly taken from Johnson~\cite{Johnson} and Pasini~\cite{Pasini-infinite}.

Let $\cP$ be a non-degenerate (ordinary) polar space. Then $\cP$ admits a Witt index if all the maximal singular subspaces have the same dimension. In view of Theorem \ref{thm21}, if a non-degenerate polar space $\cP$ admits a finite dimensional maximal singular subspace, then $\cP$ admits a Witt index.
It is worth giving a sketch of the proof of Theorem \ref{thm21}, as to make clear where the finiteness of the rank of $\cP$ plays a fundamental role. The next proposition is the main tool in the proof of Theorem \ref{thm21} and it holds for both finite rank and infinite rank polar spaces.
\begin{prop}\label{prop27}
If $M$ is a maximal singular subspace of $\cP$ and $p$ is a point of $\cP$ not in $M$, then $p^{\perp}\cap M$ is a hyperplane of $M$ and $\langle \{p\}\cup (p^{\perp}\cap M)\rangle$ is a maximal singular subspace of $\cP.$
\end{prop}
Suppose now that $\cP$ is non-degenerate and that $M$ and $M'$ are two maximal singular subspaces of $\cP.$ Put $X:=M\cap M'.$ By Proposition~\ref{prop27}, the mapping $f\colon M\rightarrow M'$ sending $x\in M$ to $x^{\perp}\cap M'$ induces an injective morphism of projective spaces from $M/X$ to the dual $(M'/X)^*$ of $M'/X$ (recall that $M$ and $M'$ are projective spaces). We have $X=\cap_{x\in M}f(x).$ Assume now that the dimension of $M$ is finite (i.e. $\di(M)<\aleph_0$). Then $X$ is the intersection of a finite number of hyperplanes of $M'$, corresponding to a basis of $M/X.$ Therefore the codimension of $X$ in $M'$ equals the codimension of $X$ in $M$ (observe that we can draw this conclusion only because $\codi_M(X)$ is finite). Hence $\di(M)=\di(M').$

If we remove the hypothesis of finiteness in Theorem \ref{thm21}, it is no more true that all maximal singular spaces are guaranteed to have the same dimension. Indeed, it is not difficult to construct polar spaces which do not admit a Witt index. See Subsection~\ref{example-infinite rank} for an example.

It is worth noting that the main achievement of the theory of polar spaces, namely the classification Theorem \ref{thm25}, holds for polar spaces of rank at least $3$ regardless of the finiteness or infiniteness of their rank.
It rather makes clear that those differences depend on the fact that sesquilinear or pseudoquadratic forms of infinite rank can behave rather differently from forms of finite rank.
\medskip

Some of the open problems related to polar spaces of infinite rank regard the validity of the following well known properties of non-degenerate polar spaces of finite rank.

\begin{itemize}
\item[(SS)] (Strong Separation Property) For every maximal singular subspace $M$, there exists a maximal singular subspace $M'$ such that $M\cap M'=\emptyset$

\item[(WS)] (Weak Separation Property) There exists at least one pair of mutually disjoint maximal singular subspaces.
\end{itemize}

Note that the strong separation properties  are proved true also for non-degenerate polar spaces of infinite rank provided there exists a maximal singular subspace of countable dimension (see Pasini~\cite{Pasini-infinite}).

\subsection{Open problems}
The main open problems regarding polar spaces of infinite rank  can be briefly summarized as the following
\begin{itemize}
\item[(O1)] Since there are examples of infinite polar spaces admitting no Witt index, are the (different) dimensions of the maximal singular subspaces in some way related? More precisely, is it possible to prove that the difference between the dimensions of the maximal polar spaces can not exceed a given infinite cardinal number?

\item[(O2)] It is an open problem to prove the validity of the strong separation proprieties for non-degenerate polar spaces of infinite not-countable rank. So far no example is known where the strong separation properties  are not true.
    \end{itemize}

\subsection{An example of polar space of infinite rank}\label{example-infinite rank}
The following construction is due to Buekenhout-Cohen~\cite{BuekCohen}. For a more detailed exposition of it I rather invite the reader to refer to~\cite{Pasini-infinite}.

 I will define a polar space of infinite rank admitting no Witt index and I will give a description of its maximal singular subspaces. Note that nearly everything that can happen in a polar space of infinite rank actually happens here, hence it can be considered as a good \rq prototype\lq\, for a polar space of infinite rank.

Let $\K$ be a given division ring, $\sigma$ an anti-automorphism of $\K$ and $\varepsilon\in \K\setminus \{0\}$ such that $\e ^{\sigma}\e=1$ and $t^{\sigma^2}=\e t\e^{-1}$ for every $t$ in $\K.$  Let $\V$ be a right $\K$-vector space of infinite dimension and let $\V^*$ be its dual, but still regarded as a right vector space over $\K$ according to the rule: $\xi \cdot t:=t^{\sigma}\xi$ for every $\xi\in \V^*$ and every $t\in \K.$ For $\xi\in \V^*$ and $x\in \V$ we put $\langle \xi,x\rangle=\xi(x)$ and we call $\langle .,.\rangle$ the {\it natural scalar product} of the pair $(\V,\V^*).$ It is not difficult to check that the form $\langle .,.\rangle$ is $\sigma$-sesquilinear.

Put $\overline{\V}:=\V\oplus\V^*$ and define a non-degenerate $(\sigma, \e)$-sesquilinear form $\Phi$ on $\overline{\V}$ as follows:

\begin{equation}
\Phi(a\oplus \alpha,b\oplus \beta):=\langle \alpha,b\rangle +\langle \beta,a\rangle^{\sigma}\e.
\end{equation}
Define a polar space $\cP$ where the points and the lines are the singular points and the singular lines of $\PG(\overline{\V})$ with respect to $\Phi.$

We immediately see that $\cP$ does not admitt any Witt index. Indeed, it is clear that $\V$ and $\V^*$  are maximal singular subspaces for $\Phi$. However, as $\di(\V)$ is infinite, $\di(\V^*)\geq 2^{\di(\V)}> \di(\V).$ So, $\Phi$ (hence $\cP$) admits no Witt index.

Denote by $\perp$ the orthogonality relation associated to $\Phi$. For a subset $A$ of $\V$ we have $A^{\perp}=\V\oplus (A^{\perp}\cap \V^*)$ and $A^{\perp}\cap \V^*$ is the subspace of $\V^*$ formed by all linear functional $\xi\in \V^*$ such that $\xi(x)=0$ for all $x\in A.$
Similarly, for $B\subseteq \V^*,$ $B^{\perp}=(B^{\perp}\cap \V) \oplus \V^*$ and $B^{\perp}\cap \V$ is formed by all vectors $x\in \V$ such that $\xi(x)=0$ for all $\xi\in B.$
Note that $A=A^{\perp \perp}$ for every subspace $A$ of $\V$.
Let now $A$ and $B$ be subspaces of $\V$ and $\V^*$ respectively such that  $A_0:=B^{\perp}\cap \V\subset A$ and $B_0:=A^{\perp}\cap \V^*\subset B.$  Then $A_0=B^{\perp}\cap A$ and $B_0=A^{\perp}\cap B.$ Hence $A_0\oplus B_0$ is singular. Moreover $\V^*/B_0$ is isomorphic to the dual $A^*$ of $A$, the quotient $\V^*/B\cong (\V^*/B_0)/(B/B_0)$ is isomorphic to $A_0^*$ and $B_0$ is isomorphic to the dual $(\V/A)^*$ of $\V/A.$ As a consequence, $B/B_0$ is isomorphic to a subspace of $(A/A_0)^*$ and $A/A_0$ is isomorphic to a subspace of $(B/B_0)^*.$ The natural scalar product $\langle ., .\rangle$ of $(\V^*,\V)$ induces a scalar product $\langle .,.\rangle_{B,A}$ for the pair $(B/B_0,A/A_0)$, defined as follows:

\[\langle \Xi, X\rangle_{B,A}=\xi(x)\,\,\,{\rm where}\,\, \xi \in \Xi\,\,{\rm and}\,\, x\in X, {\rm for} \, \Xi\in B/B_0 \,\,{\rm and}\, X\in A/A_0.\]

As $A_0\subseteq B^{\perp}$ and  $B_0\subseteq A^{\perp},$ $\xi(x)$ does not depend on the choice of $\xi\in \Xi$ and $x\in X.$ So,$ \langle .,.\rangle_{A,B}$ is well defined. Suppose that $A/A_0\cong B/B_0$  and let $f\colon A/A_0\rightarrow B/B_0$ be an isomorphism satisfying

\[\langle f(X), Y\rangle_{B,A}+\langle f(Y), X\rangle_{B,A}^{\sigma}\e=0,\,\,\,\forall X,Y\in A/A_0.\]

Put
\begin{equation}\label{max sing sub-example}
M_{A,B,f}=\{x\oplus \xi |x\in A,\xi \in f(x+A_0)\}.
\end{equation}

It can be proved that $M_{A,B,f}$ is a maximal $\Phi$-singular subspace of $\overline{\V}$ and every maximal $\Phi$-singular subspace of $\overline{\V}$ can be obtained as above, for a suitable choice of $A, B$ and $f.$ Moreover, $\di(\V)\leq \di(M_{A,B,f})\leq \di(\V^*).$



\section{Embedding polar spaces in  groups}\label{sec4}
In this section we address the problem of embeddability of polar spaces in groups. In Section \ref{sec4-1} we give the basic ingredients of a general theory of embeddings for poset geometries allowing the subgroup lattice of any group to be a feasible codomain. This abstract approach will provide a framework where both projective embeddings as defined by Ronan~\cite{Ron} and representation groups in the sense of Ivanov and Shpectorov~\cite{IS-2012} can be placed quite naturally. In Section \ref{sec4-2} we focus on embeddings of polar spaces in groups and we will report on the main achievements obtained with this respect. The source for the material of this section is Pasini~\cite{Pas}.

\subsection{Embeddings of poset-geometries}\label{sec4-1}
Let $\Gamma$ be a residually connected and firm geometry, that is the residue of any flag of corank at least $2$ is connected and every flag is contained in at least two maximal flags (see \cite{Pas-Diagram geometry}). We write $X\in \Gamma$ to say that $X$ is an element of $\Gamma.$ Given $X\in \Gamma$, we denote the type of $X$ by $t(X)$ and its residue by $Res_{\Gamma}(X)$ (or $Res(X)$). Given a subset $J\not=\emptyset$ of the set of types of $\Gamma,$ the $J$-truncation of $\Gamma$ is the geometry obtained from $\Gamma$ by removing all elements of type $j\not\in J.$
\begin{defin}
Given two geometries $\Gamma$ and $\Delta$ of rank $n$ over the same set of types and a positive integer $m<n$, an $m$-covering from $\Gamma$ to $\Delta$ is a type-preserving morphism $\varphi\colon \Gamma\rightarrow \Delta$ such that, for every flag $F$ of $\Gamma$ of corank $m$, the restriction of $\varphi$ to $Res_{\Gamma}(F)$ is an isomorphism to $Res_{\Delta}(\varphi(F)).$ Accordingly, if $\varphi\colon \Gamma'\rightarrow\Gamma$ is an $m$-covering, we say that $\Gamma'$ is an $m$-cover of $\Gamma$ and $\Gamma$ is an $m$-quotient of $\Gamma'.$

We say that a geometry $\bar{\Gamma}$ is the $m$-universal cover of $\Gamma$ if there exists an $m$-covering $\bar{\varphi}\colon \bar{\Gamma}\rightarrow \Gamma$  and for every $m$-covering $ {\varphi}'\colon  {\Gamma}'\rightarrow \Gamma$ there exists an $m$-covering $ {\psi}\colon  \bar{\Gamma}\rightarrow \Gamma '$ such that $\bar{\varphi}=\varphi ' \psi.$


If  $\Gamma$ is its own universal $m$-cover, then we say that $\Gamma$ is $m$-simply connected. The $(n-1)$-coverings are just called coverings. The universal cover of a geometry $\Gamma$ of rank $n$ is its $(n-1)$-universal cover and $\Gamma$ is said to be simply connected if it is $(n-1)$-simply connected.
 \end{defin}

A geometry $\Gamma$ is  a {\em poset-geometry} when its set of types is equipped with a total ordering $\leq$ such that, for any three elements $X,Y,Z$ of $\Gamma,$ if $t(X)\leq t(Y)\leq t(Z)$ and $Y$ is incident with both $X$ and $Z$, then $X$ is incident with $Z.$
If $X$ and $Y$ are two elements of $\Gamma,$ we write $X\leq Y$ when $X$ and $Y$ are incident and $t(X)\leq t(Y).$

I will denote by $P(\Gamma)$ the set of $0$-elements of $\Gamma$ and by $L(\Gamma)$ the set of $1$-elements of $\Gamma.$ Then $P(\Gamma)$ is called the pointset of $\Gamma$ and $L(\Gamma)$ the lineset of $\Gamma.$
If $X$ is an element of $\Gamma$, $P(X)$ denotes the set of points $p$ of $\Gamma$ such that $p\leq X$ and $L(X)$ denotes the set of lines of $\Gamma$ incident with $X.$ Moreover, we assume the following weaker version of the Intersection Property (see Pasini~\cite[Chapter 6]{Pas-Diagram geometry}): $X\leq Y$ if and only if $P(X)\subseteq P(Y).$

\begin{defin}\label{def emb}
Let $\Gamma$ be a poset geometry and $G$ a group.
An embedding of $\Gamma$ in $G$ is an injective mapping $\varepsilon \colon \Gamma\rightarrow G$   from the set of elements of $\Gamma$ to the set of proper non-trivial subgroups of $G$ such that
 \item[(E1)] for $X,Y\in \Gamma,$ we have $\varepsilon(X)\leq \varepsilon(Y)$ if and only if $X\leq Y;$
\item[(E2)] $\varepsilon (X)=\langle \varepsilon(p)\rangle_{p\in P(X)}$ for every $X$ in $\Gamma$;
\item[(E3)] $G=\langle \varepsilon(p)\rangle_{p\in P(\Gamma)}.$
\end{defin}

The group $G$ is the codomain $cod(\varepsilon)$ of $\varepsilon.$

 \begin{defin}\label{def expansion}
 Given an embedding $\varepsilon \colon \Gamma\rightarrow G$ of a poset geometry $\Gamma$ of rank $n$ in a group $G$, we define a poset-geometry $Exp(\varepsilon)$ of rank $n+1$ as follows. The points of $Exp(\varepsilon)$ are the elements of $G$ and, for $i=1,2,\dots, n$ the $i$-elements of $Exp(\varepsilon)$ are the right cosets $g\cdot \varepsilon(X),$ for $g\in G$ and $X\in \Gamma$ with $t(X)=i-1.$ The incidence relation is the inclusion between cosets and between elements and cosets. We call $Exp(\varepsilon)$ the {\em expansion} of $\Gamma$ to $G$ via $\varepsilon$.
\end{defin}

 The geometry $Exp(\varepsilon)$ is a residually connected poset-geometry. Note also that $G$, in its action by left multiplication on itself, is turned into a subgroup of $Aut(Exp(\varepsilon)).$ 

\subsection{A few important special cases of embeddings of poset-geometries}\label{special embeddings}
When $G$ is commutative then $\varepsilon$ is {\em abelian}. In particular, when $G$ is the additive group of a vector space $\V$ defined over a given division ring $\K$ and $\varepsilon (p)$ is a linear subspace of $\V$ for every $p\in P(\Gamma),$ then we say that $\varepsilon$ is a {\it $\K$-linear embedding} of $\Gamma$ in $\V$. If this is the case we write $cod(\varepsilon)=\V$ and $\varepsilon \colon \Gamma\rightarrow \V.$
Note that if $\varepsilon$ is a linear embedding, then $\varepsilon(X)$ is a linear subspace for every element $X\in \Gamma.$

If $\varepsilon$ is a $\K$-linear embedding of a poset-geometry of rank at least $2$ such that $\di(\varepsilon (p))=1$ for all points $p\in P(\Gamma)$, $\di(\varepsilon (L))=2$ for every line $L$ of $\Gamma$ and $\varepsilon (L)=\cup_{p\in P(L)}\varepsilon(p)$ for every line $L$ of $\Gamma$ then we say that $\varepsilon$ is a {\em (full) projective embedding} and we write $\varepsilon \colon \Gamma\rightarrow \PG(V).$ 

Note that if $\varepsilon$ is a projective embedding of $\Gamma$ in $\PG(\V)$ then the elements of the geometry $Exp(\varepsilon)$ defined in Definition~\ref{def expansion} are the subspaces of the affine geometry $AG(\V)$ (i.e. cosets of subspaces of $\V$).\\

Still assuming that $\Gamma$ is a poset-geometry of rank at least $2$ and $G$ is a group of $\V$, let $\varepsilon \colon \Gamma\rightarrow G$ be an embedding of $\Gamma$ and $\K$ a division ring. Suppose that two families $\{V(p)\}_{p\in P(\Gamma)}$ and $\{V(L)\}_{L\in L(\Gamma)}$ of $\K$-vector spaces are given such that the following hold
\begin{itemize}
 \item[(LP1)] for every point $p,$ $\di(V(p))=1$ and $\varepsilon(p)$ is the additive group of $V(p);$
 \item[(LP2)] for every line $L,$ $\di(V(L))=2$ and $\varepsilon(L)$ is the additive group of $V(L);$
 \item[(LP3)] for every point $p$ and every line $L$, if $p<L$ then  $V(p)$ is a subspace of $V(L);$
\item[(LP4)] for every line $L,$ $V(p)_{p\in P(L)}$ is the family of all $1$-dimensional linear subspaces of $V(L)$
\end{itemize}
Then we say that $\varepsilon$ is a {\em (full) locally $\K$-projective embedding}.

 \subsubsection{Morphisms and abstract hulls}
Given two embeddings $\varepsilon\colon\Gamma\rightarrow G$ and $\eta\colon\Gamma\rightarrow F$, a morphism from $\varepsilon$ to $\eta$ is a (group) homomorphism $f\colon G\rightarrow F$ such that, for every $X\in \Gamma,$ the restriction of $f$ to $\varepsilon (X)$ is an isomorphism to $\eta (X).$ If furthermore $f\colon G\rightarrow F$ is an isomorphism, then we say that $f$ is an isomorphism from $\varepsilon$ to $\eta.$ If a morphism exists from $\varepsilon$ to $\eta$, then we say that $\eta$ is a homorphic image of $\varepsilon$ and that $\varepsilon$ dominates $\eta.$ If there is an isomorphism from $\varepsilon$ to $\eta$ then we say that $\varepsilon$ and $\eta$ are isomorphic and we write $\varepsilon\cong \eta.$

Given a morphism $f$ from $\varepsilon$ to $\eta$, for every $X\in \Gamma$ the homomorphism $f$ maps the right cosets of $\varepsilon(X)$ in $G$ onto right cosets of $\eta(X)$ in $F$. Accordingly, $f$ defines a morphism $Exp(f)\colon Exp(\varepsilon)\rightarrow Exp(\eta)$ which is indeed a covering.

Let $\varepsilon\colon \Gamma\rightarrow G$ be an embedding. Following Ivanov~\cite{Ivanov book}, denote by $U(\varepsilon)$ the {\em universal completion} of the amalgam $\mathcal{A}(\varepsilon):=\{\varepsilon(X)\}_{X\in \Gamma}.$  Put $\phi_{\mathcal{A}(\varepsilon)}:=\{\phi_X\}_{X\in \Gamma}$,  where $\phi_X$ is the natural embedding of the group $\varepsilon (X)$ in $U(\varepsilon).$  Then  $\phi_{\mathcal{A}}$ is the natural embedding of $\mathcal{A}(\varepsilon)$ in $U(\varepsilon).$ Put $\tilde{\varepsilon}(X):=\phi_{\mathcal{A}(\varepsilon)}(\varepsilon (X))$ for every $X$ in $\Gamma.$

Then $\tilde{\varepsilon}\colon \Gamma\rightarrow U(\varepsilon)$ embeds $\Gamma$ in $U(\varepsilon).$


The canonical projection of $U(\varepsilon)$ onto $G$ induces a morphism $\pi_{\varepsilon}\colon \tilde{\varepsilon}\rightarrow \varepsilon.$
It can be proved that for every embedding $\eta$ of $\Gamma$ and every morphism $f\colon \eta\rightarrow \varepsilon$ there exists a unique morphism $g\colon \tilde{\varepsilon}\rightarrow \eta$ such that $f\circ g= \pi_{\varepsilon}.$
In view of this, we call $\tilde{\varepsilon}$ the {\em abstract hull} of $\varepsilon$ and $\pi_{\varepsilon}$  the {\em canonical projection of $\tilde{\varepsilon}$ onto $\varepsilon$}.

\begin{theo}
The geometry $Exp(\tilde{\varepsilon})$ is the universal cover of $Exp(\varepsilon)$ and, regarded $U(\varepsilon)$ and $G$ as subgroups of $Aut(Exp(\tilde{\varepsilon}))$ and $Aut(Exp( \varepsilon)),$ the group $U(\varepsilon)$ is the lifting of $G$ to $Exp(\tilde{\varepsilon}).$
\end{theo}

\begin{defin}
An embedding is abstractly dominant if it is its own hull.
\end{defin}

If we restrict to consider only special classes of embeddings, like the projective embeddings or the class of embeddings of geometries in abelian groups we can  define the concept of {\it projective hull} or {\it abelian hull} of an embedding in a similar way as we have defined the concept of abstract hull before.
More precisely, if $G$ is an abelian group and  $\varepsilon\colon \Gamma\rightarrow G$ is an embedding of $\Gamma$, then the abelian hull $\varepsilon_{ab}$ of $\varepsilon$ is the embedding $\varepsilon_{ab}\colon \Gamma\rightarrow U(\varepsilon)_{ab}$ where $U(\varepsilon)_{ab}$ is a quotient of $U(\varepsilon)$ where we have required also the commutative property to hold.

 Analogously, if $G$ is the additive group of a vector space and $\varepsilon\colon \Gamma\rightarrow G$ is a linear embedding, we consider the quotient $U(\varepsilon)_{lin}$ of $U(\varepsilon)$ requiring  it to carry the structure of a linear subspace compatible with the vector structure of $V(p)$ and $V(L)$ for $p\in P(\Gamma)$ and $L\in L(\Gamma).$ The  projective hull of the given projective embedding $\varepsilon$ of $\Gamma$ in $G$ is then $\varepsilon_{lin}\colon \Gamma\rightarrow U(\varepsilon)_{lin}.$

Consequently, we say that, in the class of projective embeddings of a geometry $\Gamma$ defined over a given division ring $\K$, where only semi-linear mappings are taken as morphisms, an embedding is {\em linearly dominant} if it is its own projective hull and that a projective embedding of $\Gamma$ is {\em linearly universal} if it is the projective hull of every projective embedding of $\Gamma$.

We would make it clear that an embedding which is linearly dominant is seldom abstractly dominant; while the abstract hull of a projective embedding is always a locally projective embedding.

Note that linearly dominant projective embeddings are often called {\em relatively universal} and linearly universal projective embeddings are often called {\em absolutely universal}  (see e.g.~\cite{Shult95}).

\subsection{Hulls of embeddings of polar spaces}\label{sec4-2}
Before focusing on polar space we need a more general definition.
Let $\Gamma$ be a thick building of connected spherical type and rank at least $2$. Given a flag $F\not=\emptyset$ and a chamber $C$ of $\Gamma$, there is a unique chamber $C_F\in Res(F)$ at minimal distance from $C$, where distance is computed in the chamber system, i.e. in the graph having as vertices the chambers and two chambers are adjacent when they intersect in a panel.
We denote the distance between $C$ and $C_F$ by $d(C,C_F).$ For every non-empty flag $X$, the distance $d(X,F)$ from $X$ to $F$ is the minimal distance $d(C,F)$ from  $F$ to a chamber $C$ containing $X.$

We say that  {\it a flag $X$ is far from} $F$ if $d(X,F)$ is maximal, compatibly with the types of $F$ and $X$.
\begin{defin}
The geometry $Far_{\Gamma}(F)$ is the substructure of $\Gamma$ formed by the elements far from $F$. Two elements $X$, $Y\in Far_{\Gamma}(F)$ are incident in $Far_{\Gamma}(F)$ if and only if they are incident in $\Gamma$ and the flag $\{X,Y\}$ is far from $F.$
\end{defin}

The geometry $Far_{\Gamma}(F)$ is residually connected except for a few cases defined over $GF(2)$ (see Blok-Brouwer~\cite{BB}).

From hereon  $\Gamma$ denotes a classical polar space of rank $n\geq 2$ and $\varepsilon\colon \rightarrow \PG(\V)$ is a full projective embedding of $\Gamma.$ As we have seen in Section~\ref{basics},  $\varepsilon(\Gamma)$ is the family of linear subspaces of $\V$ totally isotropic for a non-degenerate trace-valued reflexive form or the family of singular subspaces with respect to a non-singular pseudoquadratic form.

Denote by $\Delta:=Exp(\varepsilon)$ the expansion of $\varepsilon$  and by $\widetilde{\Delta}$ the universal cover of $\Delta.$ It is proved in Cuypers-Pasini~\cite{Cuypers-Pasini} that $\widetilde{\Delta}$ is the geometry $Far_{\Pi}(p_0)$ far from a point $p_0$ of a polar space $\Pi$ of rank $n+1$ of the same type as $\Gamma.$
Explicitly, $\Delta=\widetilde{\Delta}/E_0$ where $E_0<Aut(\widetilde{\Delta})$ is induced by the stabilizer in $Aut(\Pi)$ of all points of $\Pi$ collinear with $p_0.$\\

Denote by $\tilde{\varepsilon}$ the abstract hull of $\varepsilon.$ Then $\tilde{\varepsilon}\colon \Gamma \rightarrow U({\varepsilon})$

\begin{lemma}
The universal cover  of the expansion $ Exp(\varepsilon)$ of a full projective embedding $\varepsilon$ of a classical polar space of finite rank $n\geq 2$ is the expansion $Exp(\tilde{\varepsilon})$ of the abstract hull $\tilde{\varepsilon}$ of $\varepsilon.$
  \end{lemma}

 The following theorem completely characterizes the abstract hull $\tilde{\varepsilon}$ of a projective embedding $\varepsilon$ of a classical polar space of finite rank at least $2$ (see also Sahoo-Sastry~\cite{Sahoo-Sastry} for the case of symplectic polar spaces of odd prime power).

\begin{theo}
The codomain $U({\varepsilon})$ of $\tilde{\varepsilon}$ is isomorphic to the unipotent radical of the stabilizer of $p_0$ in $Aut(\Pi).$
\end{theo}

\begin{theo}
Suppose that $\varepsilon$ is linearly dominant. The embedding $\varepsilon$ is abstractly dominant if and only if $\varepsilon(\Gamma)$ is a quadric.
\end{theo}


\section{Dual polar spaces}\label{dual polar spaces}
In this section we focus on dual polar spaces, more precisely on embeddings of dual polar spaces. We first report on some basics about projective embeddings of a dual polar space and then we provide a survey of the known results concerning the universal and the minimal polarized projective embeddings of a classical dual polar space.

The material of this section is contained in a series of papers (\cite{BCDBP08}, \cite{C-DB-P:1}, \cite{KS}, \cite{BlCoop}, \cite{CS}, \cite{BDB-P}, \cite{DBP}).


Let $\cP$ be a classical thick polar space of rank $n \geq 2$ and let $\Delta$ be its associated dual polar space, i.e., $\Delta$ is the point-line geometry whose points, respectively lines, are the maximal, respectively next-to-maximal, singular subspaces of $\cP$  under natural incidence. If $x$ and $y$ are two points of $\Delta$, then $\d(x,y)$ denotes the distance between $x$ and $y$ in the collinearity graph of $\Delta$. If $x$ is a point of $\Delta$ and if $k \in \N$, then $\Delta_k(x)$ denotes the set of points at distance $k$ from $x$ and $x^\perp$ denotes the set of points equal to or collinear with $x$. If $X$ and $Y$ are nonempty sets of points, then $\d(X,Y)$ denotes the minimal distance between a point of $X$ and a point of $Y$. A nonempty set $X$ of points of $\Delta$ is called a {\em subspace} if every line containing two points of $X$ has all of its points in $X$. A subspace $X$ is called {\em convex} if every point on a shortest path (in the collinearity graph) between two points of $X$ is also contained in $X$.

For every point $x$ of $\Delta$, let $H_x$ denote the set of points of $\Delta$ at non-maximal distance from $x$. Then $H_x$ is a hyperplane of $\Delta$ (see~\cite{Sh-Ya}), i.e. a proper subspace of $\Delta$ meeting each line non-trivially, called a singular hyperplane of $\Delta$. It is well-known that $H_x$ is a maximal subspace of $\Delta$, see e.g. \cite{Br-Wi}.

The convex subspaces of $\Delta$ of diameter $0$ and $1$ are precisely the points and lines of $\Delta$. Convex subspaces of diameter $2$ are called {\em quads}. If $x$ is a point and $S$ is a convex subspace of $\Delta$, then $\pi_S(x)$ denotes the unique point of $S$ nearest to $x$. The point $\pi_S(x)$ is called the {\em projection} of $x$ onto $S$.

We will denote a dual polar space by putting a ``D'' in front of the name of the corresponding classical polar space. So, $\DQ(2n,\F)$, $\DQ^-(2n+1,\F)$, $\dH(2n-1,\F)$, respectively $\DW(2n-1,\F)$, denote the dual polar spaces associated with a non-singular quadric in $\PG(2n,\F)$, a non-singular elliptic quadric in $\PG(2n+1,\F)$, a non-singular hermitian variety in $\PG(2n-1,\F)$, respectively a symplectic variety in $\PG(2n-1,\F)$.

We will here consider  only embeddings in a finite-dimensional projective space. 

In Section~\ref{sec4} we already gave the definition of an embedding for an arbitrary poset-geometry. According to the terminology of Section~\ref{special embeddings}, here we will be concerned with projective embeddings of $\Delta$. We repeat the definition for more clearness.

\begin{defin}
Let $\Delta$ be a thick dual polar space of rank at least $2$. A {\em projective embedding} of $\Delta$ in a projective space $\Sigma = \PG(V)$ is an injective mapping $\varepsilon$ from the point-set $P$ of $\Delta$ to the point-set of $\Sigma$ such that:
\begin{enumerate}
\item[(PE1)] the image $\varepsilon(P)$ of $\varepsilon$ spans $\Sigma$;
\item[(PE2)] every line of $\Delta$ is mapped by $\varepsilon$ into a line of $\Sigma$;
\item[(PE3)] no two lines of $\Delta$ are mapped by $\varepsilon$ into the same line of $\Sigma$.
\end{enumerate}
\end{defin}

We will say that $\varepsilon$ is an {\em $\F$-embedding} if $\F$ is the underlying division ring of the vector space $V$. The dimensions $\dim(V)$ and $\dim(\Sigma)=\dim(V)-1$ are called the {\em vector} and {\em projective dimension} of $\varepsilon$, respectively. Note that (PE2) only says that the image $\varepsilon(L)$ of a line $L$ of $\Delta$ is contained in a line of $\Sigma$. If $\varepsilon(L)$ is a line of $\Sigma$ for every line $L$ of $\Delta$, then the embedding $\varepsilon$ is said to be {\em full}.

Since for every point $p$ of $\Delta$, the singular hyperplane $H_p$ is a maximal subspace of $\Delta$, $\varepsilon(H_p)$ spans either a hyperplane of $\Sigma$ or the whole of $\Sigma$. Following Thas and Van Maldeghem \cite{TVM}, we give the definition of polarized embedding

\begin{defin}
A projective embedding $\varepsilon$ of a dual polar space $\Delta$ in a projective space $\Sigma$ is a {\em polarized projective embedding} if $\langle \varepsilon(H_p) \rangle$ is a hyperplane of $\Sigma$ for every point $p$ of $\Delta$.
\end{defin}
If $\varepsilon$ is a polarized embedding, then $\langle \varepsilon(H_p) \rangle \cap \varepsilon(P) = \varepsilon(H_p).$ As noticed in \cite[Remarks 2 and 3]{DBP}, if the rank of $\Delta$ is sufficiently large (in any case, $n > 2$), then
$\Delta$ admits non-polarized full embeddings.

\medskip \noindent Given a projective embedding $\varepsilon:\Delta\rightarrow\Sigma$, suppose that $\alpha$ is a subspace of $\Sigma$ satisfying the following properties:
\begin{itemize}
\item[(P1)] $\alpha \cap \varepsilon(P) = \emptyset$;
\item[(P2)] $\langle \alpha,\varepsilon(x) \rangle \not= \langle \alpha,\varepsilon(y)  \rangle$ for every two distinct points $x$ and $y$ of $\Delta$.
\end{itemize}
For every point $x$ of $\Delta$, we define $\varepsilon_\alpha(x) := \langle \alpha,\varepsilon(x) \rangle$. Then $\varepsilon_\alpha$ is an embedding of $\Delta$ in $\Sigma/\alpha$, where $\Sigma/\alpha$ is formed by the subspaces of $\Sigma$ properly containing $\alpha$ (the points of $\Sigma/\alpha$ are the subspaces of $\Sigma$ containing $\alpha$ as a hyperplane). We call $\varepsilon_\alpha$ a {\em projection} or {\em quotient}  of $\varepsilon$. 
It is not difficult to see that if $\varepsilon$ is a full embedding, then the quotient $\varepsilon_\alpha$ is also full and that if $\varepsilon_\alpha$ is polarized, then $\varepsilon$ is polarized as well.

Recall from Section~\ref{sec4} that two embeddings $\varepsilon_1:\Delta \rightarrow \Sigma_1$ and $\varepsilon_2:\Delta \rightarrow \Sigma_2$ are {\em isomorphic} ($\varepsilon_1 \cong \varepsilon_2$) if there exists an isomorphism $\phi$ from $\Sigma_1$ to $\Sigma_2$ such that $\varepsilon_2(x) = \phi \circ \varepsilon_1(x)$ for every point $x$ of $\Delta$.

\medskip \noindent Suppose $\widetilde{V}$ is a vector space over the division ring $\F$. Following Cooperstein and Shult \cite{CS} we say that a full projective $\F$-embedding $\tilde{\varepsilon}:\Delta \to \widetilde{\Sigma} = \PG(\widetilde{V})$ is {\em absolutely universal projective} if for every full projective $\F$-embedding $\varepsilon$ of $\Delta$, $\varepsilon$ is isomorphic to a quotient $\tilde{\varepsilon}_{\alpha}$ of $\tilde{\varepsilon}$ for a suitable subspace $\alpha$ of $\widetilde{\Sigma}.$

The absolute universal projective embedding $\tilde{\varepsilon}$, if it exists, is uniquely determined up to isomorphisms and it is  precisely 
the linear hull of every full projective $\F$-embedding of $\Delta$ (Ronan~\cite{Ron}). Sufficient conditions for a point-line geometry to admit the full absolute universal projective  $\F$-embedding have been obtained by Kasikova and Shult \cite{KS}.

If $\varepsilon':Q \to \Sigma'$ is a full projective embedding of a thick generalized quadrangle, then by Tits \cite[8.6]{Tits}, the underlying division ring of $\Sigma'$ is completely determined by $Q$. Hence, if $\varepsilon:\Delta \to \Sigma$ is a full projective embedding of a thick dual polar space of rank $n \geq 2$, then the underlying division ring of $\Sigma$ is completely determined by $\Delta$ (since $\varepsilon$ induces a full projective embedding of each of its quads). This allows us to talk about full projective embeddings and absolutely universal embeddings of thick dual polar spaces, without mentioning the underlying division rings.


\begin{defin}
Let $\Delta$ be a dual polar space and $\varepsilon\colon \Delta \rightarrow \Sigma$ be a projective embedding of $\Delta$ in a projective space $\Sigma.$
An automorphism $g$ of $\Delta$ {\em lifts to $\Sigma$ through $\varepsilon$}
if there exists an automorphism $\varepsilon(g)$ of $\Sigma$ such that $\varepsilon(g)\varepsilon=\varepsilon g$. Clearly $\varepsilon(g)$, if it exists, is uniquely determined by $g$. We call it the {\em lifting} of $g$ to $\Sigma$ through $\varepsilon$. Let $G$ be a subgroup of the automorphism group $\Aut(\Delta)$ of $\Delta$. If all elements of $G$ lift to $\Sigma$ through $\varepsilon$ then we say that $G$ {\em lifts} to $\Sigma$, we put
$\varepsilon(G) = \{\varepsilon(g)\}_{g\in G}$ and we call $\varepsilon(G)$ the {\em lifting} of $G$ to $\Sigma$ (also, the lifting of $G$ to
$\Aut(\Sigma)$). Notice that $\varepsilon(G)$ is indeed a group. If $G$ lifts to $\Sigma$ through $\varepsilon$ then we also say that $\varepsilon$ is $G$-{\em homogeneous}.
\end{defin}

If $\varepsilon$ is $G$-homogeneous, then a subspace $U$ of $\Sigma$ is said to be {\em stabilized} by $G$ (also, to be $G$-{\em invariant}), if it is stabilized by $\varepsilon(G)$. Clearly, if a subspace $U$ of $\Sigma$ defines a quotient of $\varepsilon$ and it is $G$-invariant, then the quotient $\varepsilon_U$ is $G$-homogeneous.

We denote by $\Aut(\Delta)_0$ the normal subgroup of $\Aut(\Delta)$ generated by the root groups of the building associated to $\Delta.$ (The group $\Aut(\Delta)_0$ is in fact the largest normal simple subgroup of $\Aut(\Delta)$).
We say that a projective embedding $\varepsilon\colon \Delta \rightarrow \PG(V)$ of $\Delta$ is $\Aut(\Delta)_0$-{\em homogeneous} if $\Aut(\Delta)_0$ lifts through $\varepsilon$ to a subgroup of the full automorphism group $\PGGL(V)$ of $\PG(V)$.

The next theorems summarize the major achievements obtained for homogeneous, polarized and absolute universal embedding of dual polar spaces. In Section~\ref{cdps-embedding} we will investigate in details the situation for each classical dual polar space.

The following theorem is proved in \cite{BCDBP08}.
\begin{theo}\label{main}
If $\varepsilon$ is $\Aut(\Delta)_0$-homogeneous, then it is polarized.
\end{theo}

The next theorem follows from Buekenhout and Lef\`{e}vre \cite{BL}, Dienst \cite{Dienst} in the case of generalized quadrangles and from Kasikova and Shult \cite[Section 4.6]{KS} in the case of thick dual polar spaces of rank at least 3.

\begin{theo}\label{universalemb}
Every thick dual polar space of rank $n\geq 2$ admits the absolutely universal projective embedding, provided that it admits at least one full projective embedding.
\end{theo}

Let $\Delta$ be a thick dual polar space of rank $n \geq 2$ admitting a full polarized embedding and let $P$ denote the point-set of $\Delta$.

\begin{theo} \label{theo1-embedding}
Up to isomorphisms, there exists a unique polarized embedding $\bar{\varepsilon}$ such that every polarized embedding $\varepsilon$ of $\Delta$ has a quotient isomorphic to $\bar{\varepsilon}$. The embedding $\bar \varepsilon$ is called the {\em minimal polarized embedding} of $\Delta$.
\end{theo}

\begin{theo} \label{min-embedding}
The dimension of the minimal polarized embedding of a dual polar space of rank $n$ is at least $2^n.$
\end{theo}

\medskip \noindent Let $\tilde{\varepsilon}:\Delta\rightarrow \PG(\widetilde{V})$ denote the absolutely universal projective embedding of $\Delta$ (which exists by Theorem \ref{universalemb}) and let $\bar \varepsilon$ denote the minimal polarized embedding  (which exists by Theorem \ref{min-embedding}). Therefore, with every polarized embedding $\varepsilon$ of $\Delta$, there are associated two projections:

\[\tilde{\varepsilon} \rightarrow \varepsilon \rightarrow \bar \varepsilon.\]

 As $\tilde{\varepsilon}$ is absolutely universal, the full automorphism group $\Aut(\Delta)$ of $\Delta$ lifts through $\tilde{\varepsilon}$ to a subgroup of $\mathrm{P}\Gamma\mathrm{L}(\widetilde{V})$. Therefore, every subgroup $G$ of $\Aut(\Delta)$ lifts as well to a subgroup $\tilde{\varepsilon}(G)$ of $\mathrm{P}\Gamma\mathrm{L}(\widetilde{V})$. In particular, $\tilde{\varepsilon}$ is also $\Aut(\Delta)_0$-homogeneous. Hence it is polarized by Theorem \ref{main} (see also Cardinali, De Bruyn and Pasini \cite[Corollary 1.8]{C-DB-P:2} for a different, more straightforward proof of this claim).

Let us denote the point-set of $\Delta$ by $P$. We call $R = \bigcap_{x \in \cP} \langle \tilde{\varepsilon}(H_x) \rangle$ the {\em nucleus} of $\tilde{\varepsilon}$. By Cardinali, De Bruyn and Pasini~\cite{C-DB-P:1}, $R$ defines a quotient $\tilde{\varepsilon}_R$ of $\tilde{\varepsilon}$ which is polarized. Moreover, a projective embedding of $\Delta$ is polarized if and only if it admits a (possibly improper) quotient isomorphic to $\tilde{\varepsilon}_R$. For this reason, $\tilde{\varepsilon}_R$ is the  minimal polarized embedding of $\Delta$.
Note that $R$ might be trivial. If that is the case then $\tilde{\varepsilon}$ is the unique polarized embedding of $\Delta$.

Clearly $R$ is an $\Aut(\Delta)$-invariant subspace of $\PG(\widetilde{V})$, where we say that a subspace $U$ of $\PG(\widetilde{V})$ is $G$-{\em invariant} for a subgroup $G\leq\Aut(\Delta)$ if it is stabilized by $\tilde{\varepsilon}(G)$. As $R$ is $\Aut(\Delta)$-invariant, it is $G$-invariant for every subgroup $G$ of $\Aut(\Delta)$. In particular, $R$ is $\Aut(\Delta)_0$-invariant. Hence $\tilde{\varepsilon}_R$ is $\Aut(\Delta)_0$-homogeneous.

\begin{theo}\label{max-submod}
All $\Aut(\Delta)_0$-invariant proper subspaces of $\PG(\widetilde{V})$ are contained in $R$.
\end{theo}
In other words $R$, regarded as a subspace of $\widetilde{V}$, is the largest proper $Z$$\cdot$$\tilde{\varepsilon}(G)$-submodule of $\widetilde{V}$, where $G = \Aut(\Gamma)_0$, $Z$ stands for the center of $\mathrm{GL}(\widetilde{V})$ and $Z$$\cdot$$\tilde{\varepsilon}(G)$ is the preimage of $\tilde{\varepsilon}(G)$ by the projection of $\mathrm{GL}(\widetilde{V})$ onto $\PGL(\widetilde{V})$.

\subsection{Embeddings of classical thick dual polar spaces }\label{cdps-embedding}
\subsubsection{Embeddings  of $\DW(2n-1,\F)$}\label{symplectic dps}
Let $\W(2n-1,\F)$  be a symplectic polar space of rank $n\geq 2$ arising from a non-singular alternating form of Witt index $n$  over $\F$. Let $\Delta = \DW(2n-1,\F)$ be its dual. Then, $\Aut(\Delta)_0 = \mathrm{PSp}(2n,\F)$.

The dual polar space $\Delta$ admits the so-called {\em grassmann embedding} $\varepsilon_\gr:\Delta\rightarrow \PG(W_n)$, where $W_n=V(N,\F)\subset \bigwedge^n V $, $N = {{2n}\choose n} - {{2n}\choose{n-2}}$, mapping every totally isotropic subspace $\langle v_1,v_2,\dots, v_n\rangle$ to the point $\langle v_1\wedge v_2\wedge \dots \wedge v_n\rangle.$

Clearly, the grassmann embedding is polarized (Cardinali, De Bruyn and Pasini \cite{C-DB-P:1}) and it is proved to be universal when $\F$ is an arbitrary field different from $\F_2$ (see Cooperstein \cite{C2} for the finite case and De Bruyn and Pasini \cite{BDB-P} for the infinite case). When $\F = \F_2$, the dimension of the universal embedding of $\Delta$ is $\frac{(2^n+1)(2^{n-1}+1)}{3}$ (Li \cite{L}, Blokhuis-Brouwer \cite{BlokhuisBr}).

As far as the minimal polarized embedding is concerned, if $\char(\F) = 0$ then the grassmann embedding is the minimal polarized embedding of $\Delta.$ Indeed,
let $R$ be the nucleus of $e_\gr$, regarded as a subspace of $V = V(N,\F)$. If $\char(\F) = 0$ then $R = 0$ (see e.g. De Bruyn \cite{BartLast}).

 When $\char(\F) \neq 0$ then in general $R \neq 0$.

Suppose  $\char(\F) = 2$. By Blok, Cardinali and De Bruyn \cite{BCDB} (see also Cardinali and Lunardon \cite{CL} for the rank 3 case), if $\char(\F) = 2$ then $R$ has vector dimension $\dim(R) = N-2^n$. In this case, the minimal embedding  $e_\gr/R$  has dimension $2^n$. In particular, if $\F$ is perfect then $\Delta\cong \DQ(2n,\F)$ and $e_\gr/R$ is just the spin embedding of $\DQ(2n,\F)$.

If $\char(\F)>2$ the dimension of the minimal full polarized embedding can be computed by a recursive formula due to Premet and Suprunenko~\cite{PS1983} which indeed describes the dimensions of the modules in the composition series of $R$, where $R$ is regarded as a module for $\mathrm{PSp}(2n,\F).$ A complete geometric insight of that composition series is still missing. Some progresses aimed to explain the situation in more geometric terms have been made in Cardinali-Pasini~\cite{C-P:uniserial}.

\subsubsection{Embeddings  of $\DQ(2n,\F)$}
Let $\Q(2n,\F)$ be an orthogonal polar space of rank $n \geq 2$ of parabolic type arising from a non-singular quadratic form of Witt index $n$ over a field $\F$ and let $\Delta = \DQ(2n,\F)$ be its dual. We have $\Aut(\Delta)_0 = \mathrm{P}\Omega(2n+1,\F)$. The dual polar space $\Delta$ admits a polarized embedding $\varepsilon_\spin:\Delta\rightarrow \PG(2^n-1,\F)$, called the {\em spin embedding} of $\Delta$ (Buekenhout and Cameron \cite[Section 7]{BC}; see also De Bruyn \cite{DB:5}).

When $\char(\F)\neq 2$ the embedding $\varepsilon_\spin$ is universal (Wells \cite{Wells}) and minimal polarized. Whence it is the unique polarized embedding of $\Delta.$ The universality of the spin embedding also follows from the fact that $\Delta$ admits a generating set of size $2^n$ (Blok and Brouwer \cite{BlBr}, Cooperstein and Shult \cite{CS}).


Let now $\mathrm{char}(\F) = 2$. Suppose furthermore that $\F$ is perfect. Then $Q(2n,\F)\cong W(2n-1,\F)$, where $W(2n-1,\F)$ is the symplectic polar space of rank $n\geq 2$ arising from a non-singular alternating form of Witt index $n$. Indeed, denoted by $N_0$ the radical of the bilinear form associated to the quadratic form $q$ defining $Q(2n,\F)$, the projection from $V:=V(2n+1,\F)$ to $V/N_0 \cong V(2n-1,\F)$ induces an isomorphism from $Q(2n,\F)$ to $W(2n-1,\F)$, which can be regarded as an isomorphism from $\DQ(2n,\F)$ to $\DW(2n-1,\F)$.
We refer to Subsection~\ref{symplectic dps} for a treatment of this case.


We only recall that if $\F$ is an arbitrary field different from $\F_2$ then $\varepsilon^{sp}_n$ is universal (Cooperstein~\cite{C1} and De Bruyn-Pasini~\cite{BDB-P}).



Regarding the minimal embedding, if $\F$ is perfect then the spin embedding is the minimal polarized embedding. If $\F$ is non-perfect the minimal polarized embedding still has dimension $2^n,$ as when $\F$ is perfect.


\subsubsection{Embeddings  of $\DQ^-(2n+1,\F)$}
Let $\Q^-(2n+1,\F)$ be an orthogonal polar space of elliptic type of rank $n \geq 2$ arising from a non-singular quadratic form of Witt index $n$ over a field $\F$, which becomes a quadratic form of Witt index $n+1$ when regarded over a quadratic Galois extension $\F_1$ of $\F$. Let $\Delta = \DQ^-(2n+1,\F)$ be its dual. Then $\Aut(\Delta)_0 = \mathrm{P}\Omega^-(2n+2,\F)$. The dual polar space $\Delta$ admits a polarized embedding in $\PG(2^n-1,\F_1)$, often called the {\em spin embedding} of $\Delta$. It arises from the half-spin embedding of the half-spin geometry of $\Q^+(2n+1,\F_1)$. We refer to Cooperstein and Shult \cite{CS} and De Bruyn \cite{DB:6} for details. We shall denote this embedding by the symbol $\varepsilon^-_\spin$. The embedding $\varepsilon^-_\spin$ is universal and minimal. Indeed, $\Delta$ admits a generating set of size $2^n$ (see Cooperstein and Shult \cite{CS} for the finite case, De Bruyn \cite{DB:6} for the general case).

\subsubsection{Embeddings  of $\dH(2n-1,\F_0^2)$}
Let ${\mathrm{H}}(2n-1,\F_0^2)$ be an hermitian polar space of rank $n \geq 2$  arising from a non-singular hermitian variety of Witt index $n \geq 2$ in $\mathrm{PG}(2n-1,\F_0^2)$. Here, $\F_0$ is any field admitting a quadratic extension and $\F = \F_0^2$ is a quadratic extension of $\F_0$. 
Let $\Delta = \dH(2n-1,\F_0^2)$ be the dual of $\mathrm{H}(2n-1,\F_0^2)$. We have $\Aut(\Delta)_0 = \mathrm{PSU}(2n,\F^2_0)$. The dual polar space $\Delta$ admits a polarized embedding in $\PG(N-1,\F_0)$, where $N = {{2n}\choose n}$. (See Cooperstein \cite{C1}; also Cardinali, De Bruyn and Pasini \cite{C-DB-P:1}, De Bruyn \cite{DB:4}.) We call this embedding the {\em grassmann embedding} of $\Delta$ and we denote it by $\varepsilon^H_\gr$. 
The attribute `grassmann' is motivated by the fact that $\varepsilon^H_\gr$ arises from the usual embedding of the grassmannian of $n$-subspaces of $V(2n,\F)$ in $\PG(N-1,\F)$, via the choice of a suitable Baer subgeometry $\PG(N-1,\F_0)$ of $\PG(N-1,\F)$. If $|\F_0| > 2$ then $e^H_\gr$ is universal. Indeed, in this case $\Delta$ admits a generating set of size ${{2n}\choose n}$ (see Cooperstein \cite{C1} for the finite case and De Bruyn and Pasini \cite{BDB-P} for the general case). If $\F=\F_2$, $e^H_\gr$ is not universal. The universal embedding has dimension $\frac{(4^n+2)}{3}$ (Li~\cite{Li}).

For any choice of $\F$, the grassmann embedding $\varepsilon^H_\gr$ is the minimal embedding of $\Delta.$




\bigskip

\noindent
Ilaria Cardinali \\
Department of Information Engineering and Mathematics\\
University of Siena\\
Via Roma 56, 53100 Siena\\
ilaria.cardinali@unisi.it

\end{document}